\documentclass{article}

\input{xypic.tex} \xyoption{all}
\usepackage{amsthm,amssymb,latexsym,amsmath}

\newcommand{\Z}{\mathbb{Z}} 
\newcommand{\p}[1]{{\mathbb{P}^{#1}}} \newcommand{\pn}{{\mathbb{P}^n}}
\newcommand{\opn}{{\cal O}_{\mathbb{P}^n}} \newcommand{\ox}{{\cal O}_{X}}
\newcommand{\ov}{{\cal O}_{V}}

\newtheorem{theorem}{Theorem}

\newtheorem{proposition}[theorem]{Proposition}
\newtheorem{lemma}[theorem]{Lemma}

\newtheorem*{conjecture}{{\bf Conjecture}}
\newtheorem*{claim}{{\bf Claim}}

\begin{document}

\title{Stable bundles on 3-fold hypersurfaces}
\author{Marcos Jardim \\ IMECC - UNICAMP \\
Departamento de Matem\'atica \\ Caixa Postal 6065 \\
13083-970 Campinas-SP, Brazil}

\maketitle

\begin{abstract}
Using monads, we construct a large class of stable bundles of rank
$2$ and $3$ on 3-fold hypersurfaces, and study the set of all
possible Chern classes of stable vector bundles.

\end{abstract}


\section{Introduction}

Perhaps the most popular method of constructing rank 2 bundles over a
3-dimensional projective variety $X$ is the so-called Serre construction.
Given a local complete intersection, Cohen-Macaulay curve $C\subset X$,
let ${\cal I}_C$ denote its ideal sheaf. Then consider the extension
$$ 0 \to \ox \to E \to {\cal I}_C(k) \to 0 ~~. $$
Under some conditions on $C$, the rank 2 sheaf $E$ is locally-free; moreover,
$C$ is the zero-scheme of a section in $H^0(E)$; see \cite{Ha} for a
detailed description. In some sense, every rank 2 bundle can be obtained in
this way.

In this letter, we explore a different technique: monads. We were
motivated by a preliminary version of a paper by Douglas, Reinbacher and Yau,
who proposed, based on physical grounds, the following stronger version of the
Bogomolov inequality \cite[Conjecture 2.1]{DRY}:

\begin{conjecture}\label{dry}
Let $X$ be a non-singular, simply-connected, compact K\"ahler manifold of dimension $n$,
with K\"ahler class $H$. Assume that $X$ has trivial or ample canonical bundle. If $E$
is a $H$-stable holomorphic vector bundle over $X$ of rank $r\ge 2$, then its Chern
classes $c_1(E)$ and $c_2(E)$ satisfy the following inequality:
\begin{equation} \label{SBI}
\Delta(E) = \frac{1}{r^2}\left( 2r c_2(E) - (r-1) c_1(E)^2\right)
\cdot H^{n-2} \ge \frac{1}{12} c_2(TX) \cdot H^{n-2}
\end{equation}
\end{conjecture}

We show that this conjecture cannot be true by providing examples of stable bundles
of rank $2$ and $3$ that do not satisfy (\ref{SBI}) on hypersurfaces of degree $4$,
$5$ and $6$ within $\p4$. These counter-examples are obtained as special cases of a
more general construction of stable rank $2$ and $3$ bundles over 3-fold hypersurfaces, see
Theorem \ref{rk2} below.

This conjecture was withdrawn in a revised version of the preprint, and the counter-examples
here presented do not bear directly on the truth or falsity of the other conjectures in the revised version of \cite{DRY}. These interesting conjectures, which provide sufficient conditions for the existence of stable bundles with given Chern classes, still stand.

\paragraph{Acknowledgment.}
The author is partially supported by the CNPq grant number 300991/2004-5.


\section{Hypersurfaces and monads on hypersurfaces}

Let us begin by recalling some standard facts about hypersurfaces within
complex projective spaces.

A hypersurface $X_{(d,n)}\hookrightarrow\pn$ ($n\ge4$) of degree
$d\ge1$ is the zero locus of a section $\sigma\in H^0(\opn(d))$; for
generic $\sigma$, its zero locus is non-singular. It follows from
the Lefschetz hyperplane theorem that every hypersurface is
simply-connected and has cyclic Picard group \cite{B}. It is also
easy to see that hypersurfaces are arithmetically Cohen-Macaulay,
that is $H^p({\cal O}_{X_{(d,n)}}(k))=0$ for $1\le p\le n-1$ and all
$k\in\Z$. Finally, the restriction of the K\"ahler $\tilde{H}$ class
of $\pn$ induces a K\"ahler class $H$ on $X_{(d,n)}$, which is the
ample generator of Pic$(X_{(d,n)})$. One can show that:
$$ c_1(TX_{(d,n)}) = (n+1-d)\cdot H ~~ {\rm and} $$
\begin{equation}\label{c2}
c_2(TX_{(d,n)}) = \left( d^2-(n+1)d + \frac{1}{2}n(n+1) \right)\cdot H^2 ~~.
\end{equation}

In summary, hypersurfaces within $\pn$ with $n\ge4$ (and in fact any complete intersection
variety of dimension at least 3) do satisfy all the conditions in the Conjecture.

Fixed an ample invertible sheaf $\cal L$ with $c_1({\cal L})=H$ on a
projective variety $V$ of dimension $n$, recall that the slope $\mu
(E)$ with respect to $\cal L$ of a torsion-free sheaf $E$ on
$X_{(d,n)}$ is defined as follows:
$$ \mu (E):=\frac{c_1(E)\cdot H^{n-1}}{rk(E)} ~~. $$
We say that $E$ is stable with respect to $\cal L$ if for every
coherent subsheaf $0\ne F\hookrightarrow E$ with $0<rk(F)<rk(E)$ we
have $\mu (F) < \mu (E)$. In the case at hand, stability will always
be measured in relation to the line bundle $\ox(1)$ on the
hypersurface $X_{(d,n)}$, whose first Chern class, denoted by $H$,
is the ample generator of Pic$(X_{(d,n)})$.

A {\em linear monad} on $X_{(d,n)}$ is a complex of holomorphic
bundles of the form:
\begin{equation}\label{lmon}
0 \to \ox(-1)^{\oplus a} \stackrel{\alpha}{\longrightarrow}
\ox^{\oplus b} \stackrel{\beta}{\longrightarrow} \ox(1)^{\oplus c} \to 0 ~~,
\end{equation}
which is exact on the first and last terms. In other words, $\alpha$ is injective
and $\beta$ is surjective as bundle maps, and $\beta\alpha=0$. The holomorphic bundle
$E=\ker\beta/{\rm Im}\alpha$ is called the cohomology of the monad. Note that
$$ {\rm ch}(E) = b - a\cdot{\rm ch}(\ox(-1)) - c\cdot{\rm ch}(\ox(1)) ~~. $$
In particular,
\begin{equation}\label{ch}
{\rm rk}(E) = b-a-c ~~,~~ c_1(E) = (a-c)\cdot H ~~{\rm and}~~
c_2(E) = \frac{1}{2}(a^2-2ac+c^2+a+c)\cdot H^2 ~~,
\end{equation}
where in this case $H=c_1(\ox(1))$. The left hand side of (\ref{SBI}) is given by:
\begin{equation}\label{lm-sbi}
\Delta(E) = \frac{1}{r^2}\left( 2r c_2(E) - (r-1) c_1(E)^2
\right)\cdot H^{n-2} = \frac{b(a+c)-4ac}{(b-a-c)^2} ~~.
\end{equation}

We will also be interested in the kernel bundle $K=\ker\beta$; it has the following
topological invariants:
\begin{equation}\label{K-ch}
{\rm rk}(K) = b-c ~~,~~ c_1(K) = -c\cdot H ~~{\rm and}~~ c_2(K) =
\frac{1}{2}(c^2+c)\cdot H^2 ~~.
\end{equation}
The left hand side of (\ref{SBI}) is given by:
\begin{equation}\label{K-sbi}
\Delta(K) = \frac{1}{r^2}\left( 2r c_2(K) - (r-1) c_1(K)^2
\right)\cdot H^{n-2} =  \frac{bc}{(b-c)^2} ~~.
\end{equation}

More on linear monads and their cohomology bundles can be found at
\cite{AO,J-i,JM,JMR} and the references therein. Let us just mention
a very useful existence theorem due to Fl\o ystad in the case of
projective spaces, but easily generalizable to hypersurfaces. Below,
Fl\o ystad's original result \cite[Main Theorem]{F} is adapted to
fit our needs; the proof will not be given here, since we explicitly
establish the existence of the linear monads used in this letter.

\begin{theorem}\label{flo}
Let $X_{(d,n)}$ be a non-singular hypersurface of degree $d$ within
$\pn$, $n\ge4$. There exists a linear monad on $X$ as in
(\ref{lmon}) if and only if
\begin{itemize}
\item $ b \geq a + c + n - 2 $, if $n$ is odd;
\item $ b \geq a + c + n - 1 $, if $n$ is even
\end{itemize}\end{theorem}

Our counter-examples to Conjecture \ref{dry} will be constructed as kernel and
cohomologies of linear monads over hypersurfaces. In order to establish their
stability, we will need the following result:

\begin{theorem}\label{rk2}
Let $V$ be a 3-dimensional non-singular projective variety with
${\rm Pic}(V)=\Z$, and consider the following linear monad:
\begin{equation}\label{lm-thm}
0 \to \ov(-1)^{\oplus c} \stackrel{\alpha}{\longrightarrow}
\ov^{\oplus 2+2c} \stackrel{\beta}{\longrightarrow} \ov(1)^{\oplus
c} \to 0 ~~~~ (c\ge1)
\end{equation}
\begin{enumerate}
\item the kernel $K=\ker\beta$ is a stable rank $2+c$ bundle with
$c_1(K)=-c$ and $c_2(K)=\frac{1}{2}(c^2+c)$;
\item the cohomology $E=\ker\beta/{\rm Im}\alpha$ is a stable rank $2$ bundle
with $c_1(E)=0$ and $c_2(E)=c$.
\end{enumerate}
\end{theorem}

In Section \ref{c} below we present our counter-examples, which
arise as are special cases of Theorem \ref{rk2}. The existence of
monads of the form (\ref{lm-thm}) above for $V$ being a 3-fold
hypersurface is explicitly established in Section \ref{e}. The proof
of Theorem \ref{rk2} is left to Section \ref{p}.


\section{Counter-examples}\label{c}

Following the notation in the previous section, set $X=X_{(d,4)}$
and let $\{\sigma_1,\sigma_2,\sigma_3,\sigma_4\}$ be a basis of
$H^0(\ox(1))$. Consider the following linear monad on :
\begin{equation}\label{ex}
0 \to \ox(-1) \stackrel{\alpha}{\longrightarrow}
\ox^{\oplus 4} \stackrel{\beta}{\longrightarrow} \ox(1) \to 0
\end{equation}
with maps given by:
$$ \alpha = \left(\begin{array}{c} \sigma_1 \\ \sigma_2 \\ \sigma_3 \\ \sigma_4 \end{array}\right)
~~ {\rm and} ~~ \beta = \left(\begin{array}{cccc} -\sigma_2 ~&~
\sigma_1 ~&~ -\sigma_4 ~&~ \sigma_3 \end{array}\right) ~~. $$

It is easy to see that (\ref{ex}) is indeed a linear monad. We will
show that for $d=4,5,6$, either its kernel bundle or its cohomology
bundle will provide counter-examples to Conjecture \ref{dry}.

\subsection{Sextic within $\p4$}
Let $X=X_{(6,4)}$ be a degree $6$ hypersurface within $\p4$; notice
that $\omega_X=\ox(1)$, so that $X$ has ample canonical bundle. One
easily computes that $c_2(X)=16\cdot H^2$.

By Theorem \ref{rk2}, the cohomology of the monad (\ref{ex}) is a stable rank $2$ bundle
with $c_1=0$ and $c_2=1$. One has that
$$ \Delta(E) = \frac{1}{r^2}\left( 2r c_2(E) - (r-1) c_1(E)^2 \right)\cdot H = H^{3} $$
while
$$ \frac{1}{12} c_2(TX) \cdot H = \frac{4}{3}\cdot H^{3} ~~. $$
Therefore, the strong Bogomolov inequality (\ref{SBI}) is not
satisfied.

\subsection{Quartic within $\p4$}
Let $X=X_{(4,4)}$ be a degree $4$ hypersurface within $\p4$; notice
that $\omega_X=\ox(-1)$, so that $X$ has ample anti-canonical
bundle. One easily computes that $c_2(X)=6\cdot H^2$.

By Theorem \ref{rk2}, the kernel bundle of the monad (\ref{ex}) is a stable rank $3$
bundle with $c_1=-1$ and $c_2=1$. One has that
$$ \Delta(K) = \frac{1}{r^2}\left( 2r c_2(K) - (r-1) c_1(K)^2 \right)\cdot H = \frac{4}{9}\cdot H^{3} $$
while
$$ \frac{r^2}{12} c_2(TX) \cdot H = \frac{1}{2}\cdot H^{3} ~~. $$
Therefore, the strong Bogomolov inequality (\ref{SBI}) is not
satisfied.

\subsection{Quintic within $\p4$}
Let $X=X_{(5,4)}$ be a degree $5$ hypersurface within $\p4$; notice
that $\omega_X=\ox$, so that $X$ has trivial canonical bundle. One
easily computes that $c_2(X)=10\cdot H^2$.

By Theorem \ref{rk2}, the kernel bundle of the monad (\ref{ex}) is a stable rank $3$
bundle with $c_1=-1$ and $c_2=1$. One has that
$$ \Delta(K) = \frac{1}{r^2}\left( 2r c_2(K) - (r-1) c_1(K)^2 \right)\cdot H = \frac{4}{9}\cdot H^{3} $$
while
$$ \frac{r^2}{12} c_2(TX) \cdot H = \frac{5}{6}\cdot H^{3} ~~. $$
Therefore, the strong Bogomolov inequality (\ref{SBI}) is not
satisfied.

\subsection{Is it possible to strengthen the Bogomolov inequality?}

It is actually impossible to have an inequality of the form
\begin{equation}\label{sbi2}
\Delta(E) = \frac{1}{r^2}\left( 2r c_2(E) - (r-1) c_1(E)^2\right)
\cdot H^{n-2} \ge \kappa c_2(TX) \cdot H^{n-2}
\end{equation}
where $E$ is a stable bundle and $\kappa$ some constant, if the
underlying variety is allowed to be too general.

Indeed, as it follows from Theorem \ref{rk2} and the construction of
Section \ref{e}, given a 3-fold hypersurface $X=X_{(d,4)}$, one can
always find, for each $c\ge 1$, a stable rank $2+c$ bundle $K\to X$
with $c_1(K)=-c$ and $c_2(K)=(c^2+c)/2$, so that:
$$ \Delta(K) = \frac{(2+2c)c}{(2+c)^2} ~~. $$
Notice that the minimum value for $\Delta(K)$ is $4/9$, which occurs
exactly for $c=1$. On the other hand, by formula (\ref{c2}), the
right hand side of (\ref{sbi2}) grows quadratically with the degree
$d$.

Therefore in order for an inequality of the form (\ref{sbi2}) to
hold one must somehow restrict the type of varieties allowed, e.g.
one could take only Fano and/or Calabi-Yau varieties.

\subsection{Chern classes of stable rank 2 bundles on 3-fold hypersurfaces}
The characterization of all possible cohomology classes that arise as Chern classes
of stable bundles on a given K\"ahler manifold is not only of mathematical interest,
but it is also relevant from the point of view of physics: it amounts to describing
the set of all possible charges of BPS particles in type IIa superstring theory.

The integral cohomology ring of a 3-fold hypersurface $X=X_{(d,4)}$
is simple to describe:
$$ H^*(X,\Z) = \Z[H,L,T] / (L^2=T^2=0, H^2=dL, HL=T) ~~. $$
Notice that $H^3=dT$ and $H^4=0$. Clearly, $H$ is the generator of
$H^2(X,\Z)$, $L$ is the generator of $H^4(X,\Z)$ and $T$ is the
generator of $H^6(X,\Z)$.

Now let $E$ be a rank $r$ bundle on $X$. Recall that for any rank
$r$ bundle $E$ on a variety $X$ with cyclic Picard group, there is a
uniquely determined integer $k_E$ such that $-r+1\le c_1(E(k_E))\le
0$; the twisted bundle $E_{\rm norm}=E(k_E)$ is called the {\em
normalization} of $E$. Therefore it is enough to consider the case
when $c_1(E)=k\cdot H$ for $-r+1\le k \le 0$, and study the sets
$S_{(r,k)}(X)$ consisting of all integers $\gamma\in\Z$ for with
there exists a stable rank $r$ bundle $E$ with $c_1(E)=k\cdot H$ and
$c_2(E)=\gamma\cdot L$.

In the simplest possible case, provided by $d=1$ (so that $X=\p3$)
and $r=2$, this problem was completely solved by Hartshorne in
\cite{Ha}. He proved that $S_{(2,0)}(\p3)$ consists of all positive
integers, while $S_{(2,-1)}(\p3)$ consists of all positive even
integers. As far as it is known to the author, Hartshorne's result
has not been generalized for other 3-folds.


As a consequence of Theorem \ref{rk2}, we have:

\begin{lemma}
For every positive integer $c\ge1$, $cd\in S_{(2,0)}(X_{(d,4)})$.
\end{lemma}

Based on Hartshorne's result mentioned above, it seems reasonable to conjecture that
$S_{(2,0)}(X_{(d,4)})$ consists exactly of all positive multiples of $d$.

The monad construction does not yield stable rank 2 bundles with odd
first Chern class; to construct those, one needs a variation of the
Serre construction, which provides a 1-1 correspondence between rank
2 bundles and codimension 2 subvarieties on $\p3$; see Hartshorne's
paper \cite{Ha}.


\section{Existence of linear monads on 3-fold hypersurfaces}\label{e}

Let $X=X_{(d,4)}$ be a hypersurface of degree $d$ within $\p4$; as
above let $\{\sigma_1,\sigma_2,\sigma_3,\sigma_4\}$ be a basis of
$H^0(\ox(1))$. We will now explicitly establish the existence of
linear monads of the form
$$ 0 \to \ox(-1)^{\oplus c} \stackrel{\alpha}{\longrightarrow}
\ov^{\oplus 2+2c} \stackrel{\beta}{\longrightarrow} \ox(1)^{\oplus
c} \to 0 ~~~~ (c\ge1) ~~. $$

Consider the $c\times(c+1)$ matrices:
$$ B_1 = \left( \begin{array}{ccccc}
\sigma_1 & \sigma_2 & & & \\ & \sigma_1 & \sigma_2 & & \\ & & \ddots & \ddots & \\
& & & \sigma_1 & \sigma_2
\end{array} \right) ~~~~
B_2 = \left( \begin{array}{ccccc} \sigma_3 & \sigma_4 & & & \\ &
\sigma_3 & \sigma_4 & & \\ & & \ddots & \ddots & \\ & & & \sigma_3 &
\sigma_4
\end{array} \right) ~~, $$
and the $(c+1)\times c$ matrices:
$$ A_1 = \left( \begin{array}{cccc}
\sigma_2 & & & \\ \sigma_1 & \sigma_2 & & \\ & \ddots & \ddots & \\
& & \sigma_1 & \sigma_2 \\ & & & \sigma_1
\end{array} \right) ~~~~
A_2 = \left( \begin{array}{cccc} \sigma_4 & & & \\ \sigma_3 & \sigma_4 & & \\
& \ddots & \ddots & \\ & & \sigma_3 & \sigma_4 \\ & & & \sigma_3
\end{array} \right) ~~, $$
Notice that all four matrices have maximal rank $c$.
It easy to check that:
$$ B_1A_2 = B_2A_1 = \left( \begin{array}{cccccc}
\phi_1 & \phi_2 & & & & \\ & \phi_0 & \phi_1 & \phi_2 & & \\
& & \ddots & \ddots & \ddots & \\ & & & & \phi_0 & \phi_1
\end{array} \right) ~~, $$
where $\phi_0 = \sigma_1\sigma_3$,
$\phi_1=\sigma_1\sigma_4+\sigma_2\sigma_3$ and $\phi_2 =
\sigma_2\sigma_4$.

Now form the linear monad:
$$ 0 \to \ox(-1)^{\oplus c} \stackrel{\alpha}{\to} \ox^{\oplus 2+2c}
\stackrel{\beta}{\to} \ox(1)^{\oplus c} \to 0 $$
where the maps $\alpha$ and $\beta$ are given by:
$$ \beta = \left( \begin{array}{c} B_1 \\ B_2 \end{array}\right)
~~{\rm and}~~ \alpha = \left( \begin{array}{cc} A_2 \\ -A_1
\end{array}\right) $$ Clearly, both maps are of maximal rank $c$ for
every point in $X$, and $\beta\alpha=B_1A_2-B_2A_1=0$.


\section{Proof of Theorem \ref{rk2}}\label{p}

The proof is based on a very useful criterion (due to Hoppe) to decide whether a
bundle on a variety with cyclic Picard group is stable. We set $E_{norm}:=E(k_E)$
and we call $E$ normalized if $E=E_{\rm norm}$. We then have the following criterion.

\begin{proposition}\label{hoppe} (\cite[Lemma 2.6]{H})
Let $E$ be a rank $r$ holomorphic vector bundle on a variety $X$ with Pic$(X)=\Z$.
If $H^0((\wedge ^qE)_{\rm norm})=0$ for $1\leq q\leq r-1$, then
$E$ is stable.
\end{proposition}

Our argument follows \cite[Theorem 2.8]{AO}. Consider the linear monad
$$ 0 \to \ox(-1)^{\oplus c} \stackrel{\alpha}{\to} \ox^{\oplus 2+2c}
\stackrel{\beta}{\to} \ox(1)^{\oplus c} \to 0 ~~; $$
setting $K=\ker\beta$; one has the sequences:
\begin{equation} \label{ker1}
0 \to K \to \ox^{\oplus 2+2c} \stackrel{\beta}{\longrightarrow}
\ox(1)^{\oplus c} \to 0 ~~ {\rm and}
\end{equation}
\begin{equation} \label{ker2}
0 \to \ox(-1)^{\oplus c} \stackrel{\alpha}{\longrightarrow} K \to E \to 0 ~~.
\end{equation}
First, we will show that the kernel bundle $K$ is stable. That implies that $K$ is
simple, which in turn implies that cohomology bundle $E$ is simple
Since any simple rank 2 bundle is stable, we conclude that $E$ is also stable.

Recall that one can associate to the short exact sequence of locally-free sheaves
$0\to A\to B\to C\to 0$ two long exact sequences of symmetric and exterior powers:
\begin{equation} \label{one}
0 \to \wedge^qA \to \wedge^qB \to \wedge^{q-1}A\otimes C \to \cdots \to
B\otimes S^{q-1}C \to S^{q}C \to 0
\end{equation}
\begin{equation} \label{two}
0 \to S^qA \to S^{q-1}A\otimes B \to \cdots \to
A\otimes\wedge^{q-1}B \to \wedge^qB \to \wedge^{q}C \to 0
\end{equation}
In what follows, $\mu(F)=c_1(F)/{\rm rk}(F)$ is the slope of the sheaf $F$,
as usual.

Finally, notice that $H^p(\ox(k))=0$ for $p\ge2$ and $k\ge-1$, by
the Kodaira vanishing theorem.

\begin{claim}
$K$ is stable.
\end{claim}

From the sequence dual to sequence (\ref{ker1}), we get that:
$$ \mu(K^*) = \frac{c}{c+2} ~~\Longrightarrow~~ \mu(\wedge^qK^*) = \frac{qc}{c+2} $$
so that $(\wedge ^qK^*)_{\rm norm}=\wedge^qK^*(k)$ for some
$k\le-1$, and if $H^0(\wedge ^qK^*(-1))=0$, then $H^0((\wedge
^qK^*)_{\rm norm})=0$.

The vanishing of $h^0(K^*(-1))$ (i.e $q=1$) is obvious from the dual to
sequence (\ref{ker1}). For the case $q=2$, start from the dual to (\ref{ker1})
and consider the associated sequence
$$ 0\to S^2(\ox(-1)^{\oplus c}) \to \ox(-1)^{\oplus c}\otimes\ox^{\oplus 2c+2}
\to \wedge^2(\ox^{\oplus 2c+2}) \to \wedge^2K^* \to 0 ~~. $$
Twist it by $\ox(-1)$ and break it into two short exact sequences:
$$ 0 \to \ox(-3)^{\oplus{{c+1}\choose{2}}} \to \ox(-2)^{\oplus 2c^2+2c} \to Q \to 0 $$
$$ 0 \to Q \to \ox(-1)^{\oplus{{2c+2}\choose{2}}} \to \wedge^2K^*(-1) \to 0 $$
Passing to cohomology, we get $H^0(\wedge ^2K^*(-1))=H^1(Q)=0$.

Now set $q=3+t$ for $t=0,1,\cdots,c-2$ and note that
$$ \mu(\wedge^{3+t}K^*(-t-1)) = \frac{(3+t)c}{c+2}-t-1 = 2\frac{c-t-1}{c+2}>0 ~~. $$
Thus $(\wedge^{3+t}K^*)_{\rm norm}=\wedge^{3+t}K^*(k)$ for some $k\le-t-2$, and
if $H^0(\wedge^{3+t}K^*(-t-2))=0$, then $H^0((\wedge^{3+t}K^*)_{\rm norm})=0$.

We show that $H^0(\wedge^{3+t}K^*(-t-2))=0$ by induction on $t$. From the dual
to sequence (\ref{ker2}) we get, after twisting by $\ox(-2)$:
$$ 0 \to \wedge^{3}K^*(-2) \to \wedge^{2}K^*(-1)^{\oplus c} \to \cdots $$
since $\wedge^{3}E^*=0$ because $E$ has rank 2. Passing to cohomology, we get
that $H^0(\wedge^{3}K^*(-2))=0$, since, as we have seen above, $H^0(\wedge^{2}K^*(-1))=0$.
This proves the statement for $t=0$.

By the same token, we get from the dual to sequence (\ref{ker2})
after twisting by $\ox(-2-t)$:
$$ 0 \to \wedge^{3+t}K^*(-2-t) \to \wedge^{2+t}K^*(-t-1)^{\oplus c} \to \cdots ~~.$$
Passing to cohomology, we get
$$ H^0(\wedge^{2+t}K^*(-t-1))=0 ~~ \Rightarrow ~~ H^0(\wedge^{3+t}K^*(-t-2))=0 $$
which is the induction step we needed.

In summary, we have shown that $H^0((\wedge ^qK^*)_{\rm norm})=0$ for
$1\leq q\leq c+1$, thus by (\ref{hoppe}) we complete the proof of the
claim.

\begin{claim}
$E$ is simple, hence stable.
\end{claim}

Tensoring by $E$ the sequence dual to (\ref{ker2}) we get
\begin{equation}\label{ee}
0\to H^0(E^*\otimes E)\to H^0(K^*\otimes E)\to\cdots ~~.
\end{equation}
Now tensoring (\ref{ker2}) by $K^*$ we get:
$$ H^0(K^*(-1))^{\oplus c} \to H^0(K^*\otimes K)\to
H^0(K^*\otimes E) \to H^1(K^*(-1))^{\oplus c} ~~. $$

But it follows from the dual of sequence (\ref{ker1}) twisted by
$\ox(-1)$ that $h^0(K^*(-1))=h^1(K^*(-1))=0$; thus $h^0(E^*\otimes
E)=1$ because $K$ is simple. But $E$ has rank 2, thus $E$ is stable,
as desired.

This completes the proof of the Theorem \ref{rk2}.
\hfill$\Box$


 \end{document}